\documentclass[10pt,a4paper]{article}
\usepackage[T1]{fontenc}
\usepackage{graphicx}
\usepackage{amssymb}
\usepackage{amsthm}
\usepackage{amsmath}
\usepackage{hyperref}
\usepackage{multirow}
\usepackage{subcaption}
\usepackage{array}
\usepackage{booktabs}
\usepackage{float}
\usepackage{algorithm}
\usepackage[noend]{algpseudocode}
\PassOptionsToPackage{numbers}{natbib}
\usepackage[preprint]{neurips_2024}
\newtheorem{theorem}{Theorem}[section]

\newtheorem{proposition}[theorem]{Proposition}

\author{%
	Xu Duan$^1$, Dongmei Chen$^{1}$\\
	$^1$The University of Texas at Austin\\
	\texttt{\{xu.duan, dmchen\}@utexas.edu}}
\title{Cosine-Similarity Methods for Efficient Training and Sampling in High-Dimensional Latent Spaces}
\begin{document}
	\maketitle

	\begin{abstract}
	Latent generative models are increasingly shifting from traditional VAEs toward representation autoencoders and semantically aligned latent spaces, which lift images into higher-dimensional feature domains where semantic factors become more separable. Yet these spaces also contain geometric regularities that existing methods do not fully exploit—particularly in the directional relationships between features.
	We introduce a cosine-similarity–based mechanism that improves both training and sampling by selecting couplings that produce cleaner, less entangled velocity fields. This simple alignment reduces gradient noise, accelerates convergence, and improves sample fidelity. Building on this idea, we develop cosine-similarity–based fine-tuning and time-scheduling strategies that reduce the FID of an 800-epoch RAE from 11.99 to 8.60. Furthermore, by formulating an optimal-transport coupling using a cosine cost, a single-epoch fine-tuning step at the 20-epoch checkpoint reaches 3.30 FID—matching the performance of the 80-epoch baseline.
	\end{abstract}
	
	\section{Introduction}
	Pixel-space diffusion models remain computationally demanding, motivating the shift toward compact latent-space generation. Models such as Latent Diffusion Models (LDMs) \cite{rombachHighResolutionImageSynthesis2022} and Diffusion/Flow Transformers (DiT/SiT) \cite{peeblesScalableDiffusionModels2023, maSiTExploringFlow2024} demonstrate that learning and sampling in a low-dimensional latent space yields higher visual fidelity and efficiency compared to pixel-space diffusion.
	
	Concurrently, advances in self-supervised representation learning—including DINO, DINOv2,  MAE, and MOCOV3 \cite{assranSelfSupervisedLearningImages2023, oquabDINOv2LearningRobust2024, heMaskedAutoencodersAre2021, chenEmpiricalStudyTraining2021}—produce encoders whose feature spaces exhibit robust semantic structure and strong generalization. These learned representations naturally raise the question of whether generative models should incorporate such semantics into their latent spaces.
	
	Recent work increasingly shows that generation benefits from semantically constrained latent variables. Yao et al. \cite{yaoReconstructionVsGeneration2025} highlight that generation in unconstrained VAE latents leads to semantic drift, while REPA-style methods \cite{yuRepresentationAlignmentGeneration2025, lengREPAEUnlockingVAE2025} explicitly align generative outputs with self-supervised visual features to improve consistency and reduce distortions.
	
	Other approaches adapt the generative noise itself. DNA-Edit \cite{xieDNAEditDirectNoise2025} optimizes noise vectors along semantic directions to improve realism. Tong et al. \cite{tongImprovingGeneralizingFlowbased2024} further formalize noise design as an optimal-transport problem, aligning data with structured Gaussian priors to produce smoother training dynamics. Li et al. propose Immiscible Diffusion \cite{liImmiscibleDiffusionAccelerating2024}, arguing that fully mixing all images across the noise space complicates denoising. Their assignment-then-diffusion strategy pairs each image with nearby noise before diffusion, reducing noise–data entanglement. This simple modification preserves Gaussian structure while significantly accelerating training and improving fidelity.
	
	Despite recent progress, high-dimensional representation spaces still contain rich geometric structure that generative models rarely capitalize on. We develop a cosine-similarity framework that exploits this structure by coupling training targets and sampling trajectories according to their directional agreement in feature space. This alignment not only accelerates convergence but also improves sample fidelity.
	Additionally, we introduce cosine-driven scheduling strategies—applicable during both training and ODE-based generation—that adapt step sizes based on semantic proximity. Finally, we construct an optimal-transport formulation using a cosine cost, enabling a one-epoch fine-tuning procedure that pushes a low-epoch checkpoint to near-saturated performance.
	
	Together, these techniques illustrate that carefully leveraging semantic geometry can yield substantial gains in generative modeling—without modifying architectures or increasing inference complexity.
	\begin{figure}
		\includegraphics[scale=0.6]{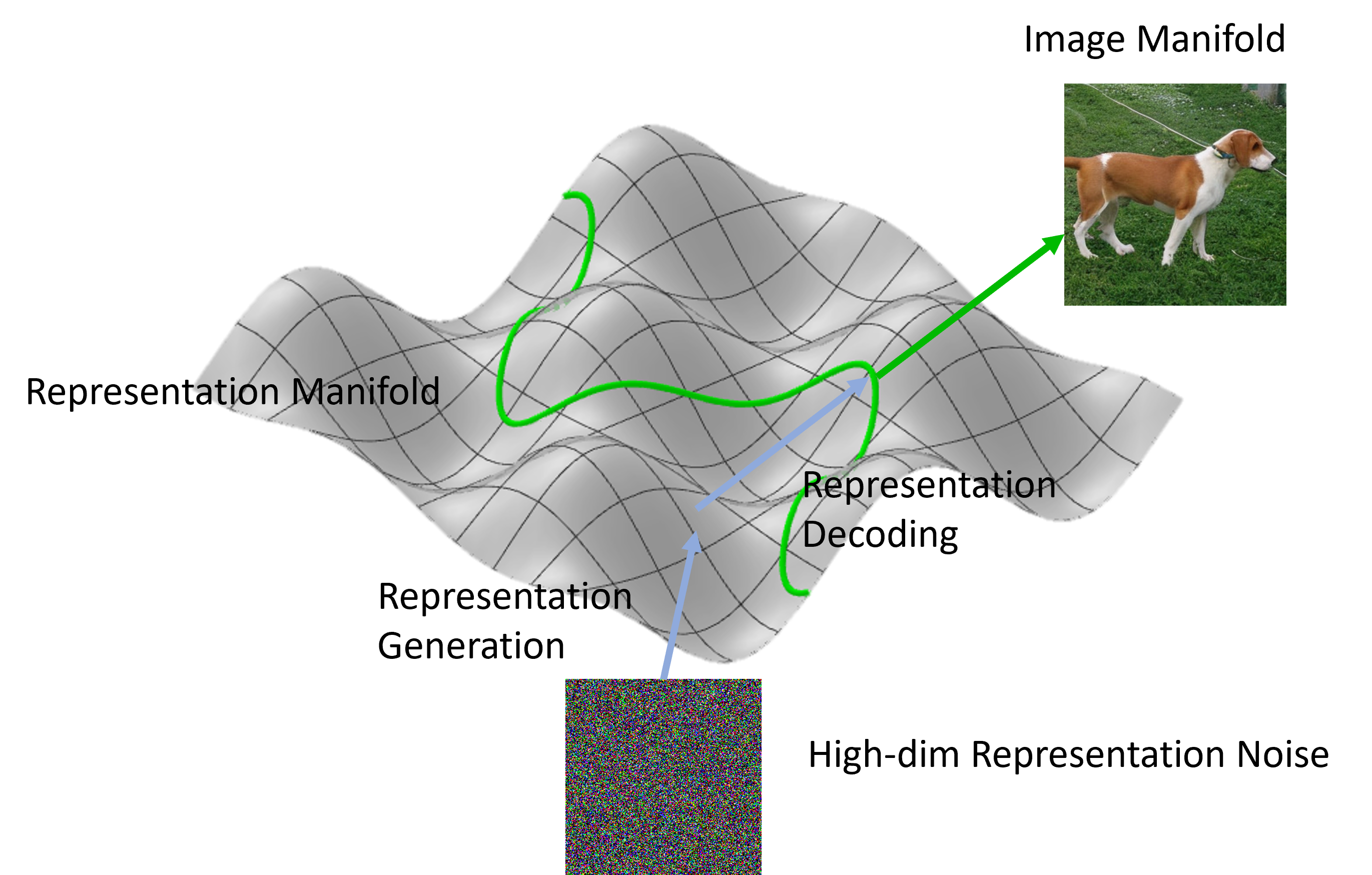}
		\caption{Scratch of High Dimension Representation Space Generation}
	\end{figure}

	\section{Related Work}
	Generative models for image generation.
	Classical diffusion models such as DDPM \cite{hoDenoisingDiffusionProbabilistic2020} and DDIM \cite{songDenoisingDiffusionImplicit2020} perform progressive noise removal over hundreds or thousands of steps, resulting in substantial computational cost during both training and inference. Latent diffusion models (LDMs) \cite{rombachHighResolutionImageSynthesis2022} mitigate this by operating in a compressed latent space using pretrained autoencoders, while transformer-based architectures—including DiT \cite{peeblesScalableDiffusionModels2023} and SiT \cite{maSiTExploringFlow2024}—enhance global context modeling via self-attention mechanisms. Despite these architectural advances, diffusion training remains slow due to the need to learn denoising scores across an entire noise schedule.
	
	A large body of work seeks to accelerate diffusion through improved sampling or revised training procedures. Sampling-side acceleration includes higher-order solvers \cite{luDPMSolverFastODE2022, luDPMSolverFastSolver2022, zhangFastSamplingDiffusion2022}, consistency distillation \cite{songConsistencyModels2023}, and one-step or few-step student models \cite{gengMeanFlowsOnestep2025, yinOnestepDiffusionDistribution2023}, though these often require expensive teacher models or can introduce instability. Training-side acceleration employs masked or partial-timestep training \cite{leiMaskedDiffusionModels2023, zhengFastTrainingDiffusion2023}, multi-resolution optimization \cite{yellapragadaZoomLDMLatentDiffusion2024}, or architectural modifications, but these typically increase implementation complexity. Other approaches compress the latent space through quantization or semantic bottlenecks \cite{xuMSFEfficientDiffusion2025}, sometimes at the cost of generation fidelity.
	
	Image generation with semantic representations.
	Recent work explores coupling generative models with pretrained representation encoders. VA-VAE \cite{yaoReconstructionVsGeneration2025} aligns VAE latents with semantic features, and MAETok \cite{chenMaskedAutoencodersAre2025}, DC-AE 1.5 \cite{chenDCAE15Accelerating2025}, and l-DeTok \cite{yangLatentDenoisingMakes2025} incorporate MAE- or DAE-style objectives \cite{vincentExtractingComposingRobust2008} into VAE training. While alignment improves both reconstruction and sample quality, the reliance on heavily compressed, low-dimensional latents limits fidelity and the richness of learned representations.
	
	In contrast, we reconstruct directly from representation-encoder features without additional compression, and we show that a simple ViT-based decoder on top of frozen features achieves reconstruction performance comparable to or exceeding SD-VAE \cite{rombachHighResolutionImageSynthesis2022} while preserving significantly stronger representations.
	
	Representation for generation.
	A parallel line of work investigates using semantic representations to improve generative modeling itself. REPA \cite{yuRepresentationAlignmentGeneration2025} accelerates DiT convergence by aligning its middle block with representation-encoder features; REPA-E \cite{lengREPAEUnlockingVAE2025} extends this idea by also unfreezing the VAE. DDT \cite{wangDDTDecoupledDiffusion2025} further improves convergence by decoupling DiT into an encoder–decoder structure and applying REPA loss to the encoder output. REG \cite{wuRepresentationEntanglementGeneration2025} introduces a learnable token into the DiT sequence and aligns it with semantic features, while ReDi \cite{kouzelisBoostingGenerativeImage2025} jointly models VAE latents and PCA-compressed DINOv2 features within a diffusion framework. Orthogonally, \cite{wangDiffuseDisperseImage2025} leverage contrastive-learning-based loss designs to improve image quality in a plug-and-play fashion without modifying the network architecture.

	\section{Background}
	Our work builds on RAE \cite{zhengDiffusionTransformersRepresentation2025}.
	We begin by introducing the relevant preliminaries.
	Flow and diffusion models both leverage stochastic processes to gradually transform Gaussian noise
	$\epsilon \sim \mathcal{N}(0, I)$ into data samples $\mathbf{x}_*$. This process can be unified as
	\begin{equation}
		\mathbf{x}_t = \alpha_t \mathbf{x}_* + \sigma_t \epsilon,
	\end{equation}
	where $\alpha_t$ is a decreasing and $\sigma_t$ an increasing function of time $t$. Flow-based models
	typically interpolate between noise and data over a finite interval, while diffusion models define a
	forward stochastic differential equation (SDE) that converges to a Gaussian distribution as
	$t \rightarrow \infty$.
	
	Sampling from these models can be achieved via either a reverse-time SDE or a probability flow
	ordinary differential equation (ODE), both of which yield the same marginal distributions for
	$\mathbf{x}_t$. The probability flow ODE is:
	\begin{equation}
		\dot{\mathbf{x}}_t = \mathbf{v}(\mathbf{x}_t, t),
	\end{equation}
	where the velocity field $\mathbf{v}(\mathbf{x}, t)$ can be formulated by the conditional
	expectation:
	\begin{equation}
		\mathbf{v}(\mathbf{x}, t)
		= \mathbb{E}[\dot{\mathbf{x}}_t \mid \mathbf{x}_t = \mathbf{x}]
		= \dot{\alpha}_t \mathbb{E}[\mathbf{x}_* \mid \mathbf{x}_t = \mathbf{x}]
		+ \dot{\sigma}_t \mathbb{E}[\epsilon \mid \mathbf{x}_t = \mathbf{x}].
	\end{equation}
	
	To synthesize data, we can integrate Eqn.~(3) in reverse time, initializing from
	$\mathbf{X}_T = \epsilon$ where $\epsilon \sim \mathcal{N}(0, I)$. This process yields samples from
	$p_0(\mathbf{x})$, serving as an approximation to the true data distribution $p(\mathbf{x})$. This
	velocity can be estimated by a model $\mathbf{v}_\theta(\mathbf{x}_t, t)$, which is trained to
	minimize the following loss function:
	\begin{equation}
		\mathcal{L}_v(\theta)
		= \int_0^T \mathbb{E}\!\left[
		\left\|
		\mathbf{v}_\theta(\mathbf{x}_t, t)
		- \dot{\alpha}_t \mathbf{x}_*
		- \dot{\sigma}_t \epsilon
		\right\|^2
		\right] dt.
	\end{equation}
	
	The reverse-time SDE can describe the probability distribution $p_t(\mathbf{x})$ of $\mathbf{x}_t$
	at time $t$, which can be expressed as:
	\begin{equation}
		d\mathbf{x}_t
		= \mathbf{v}(\mathbf{x}_t, t) dt
		- \frac{1}{2} w_t \mathbf{s}(\mathbf{x}_t, t) dt
		+ \sqrt{w_t}\, d\overline{\mathbf{W}}_t,
	\end{equation}
	with $\mathbf{s}(\mathbf{x}, t)$ denoting the score that can be computed via the conditional
	expectation:
	\begin{equation}
		\mathbf{s}(\mathbf{x}_t, t)
		= -\sigma_t^{-1} \mathbb{E}[\epsilon \mid \mathbf{x}_t = \mathbf{x}].
	\end{equation}
	The score can be reformulated in terms of the velocity $\mathbf{v}(\mathbf{x}, t)$:
	\begin{equation}
		\mathbf{s}(\mathbf{x}, t)
		= \sigma_t^{-1}
		\cdot
		\frac{
			\alpha_t \mathbf{v}(\mathbf{x}, t) - \dot{\alpha}_t \mathbf{x}
		}{
			\alpha_t \dot{\sigma}_t - \dot{\alpha}_t \sigma_t
		}.
	\end{equation}
	
	We can learn the velocity field $\mathbf{v}(\mathbf{x}, t)$ and use it to compute the score
	$\mathbf{s}(\mathbf{x}, t)$ when using an SDE for sampling.

	\section{Cosine Similarity Based Sampling Techniques}
	Esser et al.~\cite{esserScalingRectifiedFlow2024} propose a noise–to–data time–warping
	mechanism based on the signal-to-noise ratio (SNR), designed to slow down integration in
	noise-dominated regions while accelerating it near the data manifold. Their formulation
	introduces a nonlinear reparameterization of time,
	\begin{equation}
		t_m = 
		\frac{\sqrt{m/n}\; t_n}{1 + (\sqrt{m/n} - 1)\, t_n},
		\label{eq:esser-time-shift}
	\end{equation}
	which redistributes integration effort toward regions where the model dynamics exhibit
	rapid transitions.
	
	\medskip
	
	Inspired by this principle, we design an adaptive sampler that adjusts the integration
	step size using the cosine similarity between the current state $\mathbf{x}$ and the 
	instantaneous drift field $\mathbf{f}(\mathbf{x}, t)$. Intuitively, the cosine similarity 
	serves as a proxy for the ``alignment'' between position and velocity: 
	large negative similarity indicates that the drift induces a sharp directional change, 
	suggesting that smaller steps are necessary, whereas high similarity corresponds to 
	smooth regions in which larger steps are safe.
	
\begin{figure}[t]
		\centering
		\begin{subfigure}{0.48\textwidth}
			\centering
			\includegraphics[width=\linewidth]{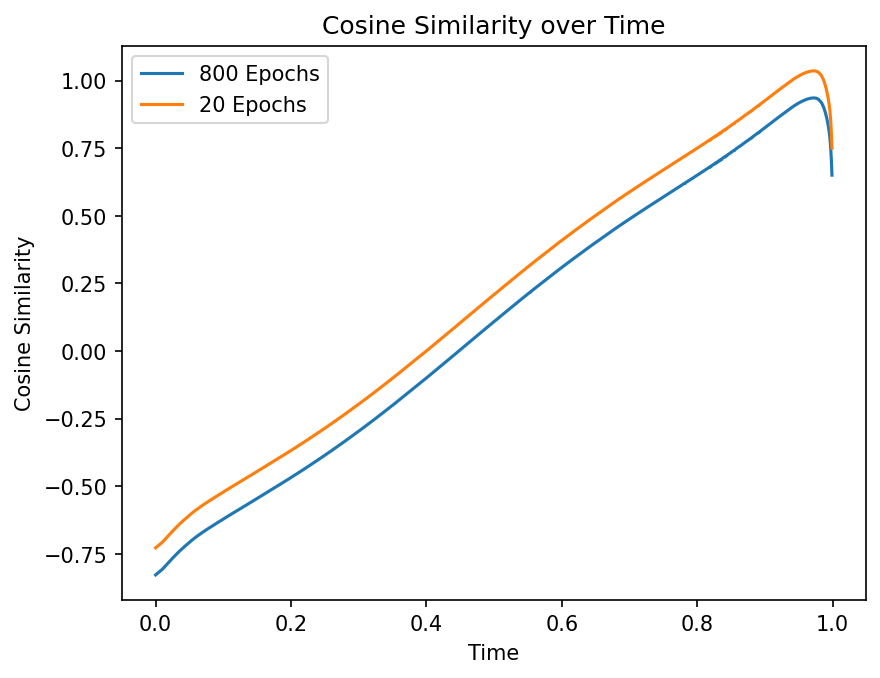}
			\caption{Cosine similarity between velocity and position over time.}
			\label{fig:cos-sim}
		\end{subfigure}
		\hfill
		\begin{subfigure}{0.48\textwidth}
			\centering
			\includegraphics[width=\linewidth]{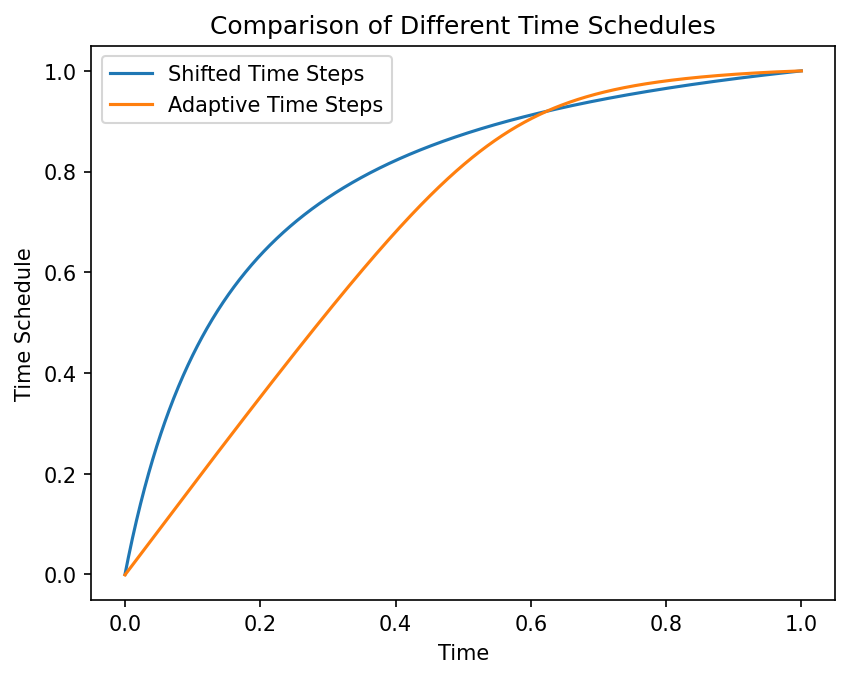}
			\caption{Different time-stepping schedules used during integration.}
			\label{fig:time-schedule}
		\end{subfigure}
		
		\caption{Comparison of cosine-similarity trends (left) and the time-step scheduling strategies (right).}
		\label{fig:cos-sim-time-schedule}
	\end{figure}

	As shown in Fig.~\ref{fig:cos-sim}, the cosine similarity between the state $\mathbf{x}_t$ and the drift field exhibits a characteristic shape: when $t$ is small (i.e., the sample is still close to real data), the similarity is negative and its magnitude is relatively large. This indicates that the drift direction is sharply misaligned with the current state, and therefore smaller integration steps are desirable. As $t$ increases and the trajectory moves deeper into the noise regime, the cosine similarity gradually increases toward~1, corresponding to smoother dynamics that allow larger step sizes.
	
	Given the cosine similarity $\mathrm{cos}$ at each step, we construct a normalized control signal
	\begin{equation}
		\tilde{c}
		= \frac{1 - \mathrm{cos}(\mathbf{x}, \mathbf{f})}{2}
		\;\in\; [0,1],
	\end{equation}
	which is larger whenever the drift direction differs sharply from the current state.  
	To obtain a sharper nonlinear response, we apply a sigmoid transformation
	\begin{equation}
		\alpha = \sigma\!\left(10(\tilde{c} - 0.5)\right),
	\end{equation}
	and define the adaptive time step
	\begin{equation}
		\Delta t
		= \Delta t_{\min}
		+ (\Delta t_{\max} - \Delta t_{\min}) \, \alpha .
	\end{equation}
	
	This scheme increases the step size in smooth regions and reduces it when the sampler encounters rapid geometric changes, providing a direct, data-adaptive alternative to SNR-based time reparameterization.

	Fig.~\ref{fig:time-schedule} compares our cosine-similarity–based time schedule with the SNR-based time shift proposed by Esser et al.~\cite{esserScalingRectifiedFlow2024}. While Esser’s schedule focuses on reallocating computation according to the global SNR profile, our method produces a more adaptive progression that directly reflects the instantaneous geometry of the learned vector field. This yields a time schedule that is more responsive in regions of rapid directional change and more permissive where the flow is smoother.

	\begin{figure}[H]
		\centering
		
		\begin{minipage}{0.45\textwidth}
      	Fig.~\ref{fig:fid-timeschedule} reports FID scores under different
		numbers of sampling steps for both the RAE time schedule and our
		cosine-similarity–based schedule.
		At 800 training epochs, our method consistently achieves lower FID
		across all step budgets, with the largest gain appearing in the
		low-step regime (e.g., $5$ steps: $8.6$ vs.\ $11.91$).
		This indicates that the adaptive schedule allocates computation more
		effectively when the sampler must operate under severe step
		constraints.
		
		A similar trend appears at 20 training epochs.  
		Although absolute performance is lower due to early-stage
		representations, our schedule still offers a uniform improvement
		(e.g., $5$ steps: $9.7$ vs.\ $10.33$).
		Overall, these results show that cosine-similarity time warping
		provides a more efficient integration trajectory than the SNR-based
		RAE schedule.
		\end{minipage}
		\hfill
		\begin{minipage}{0.52\textwidth}
			\centering
			\includegraphics[width=\linewidth]{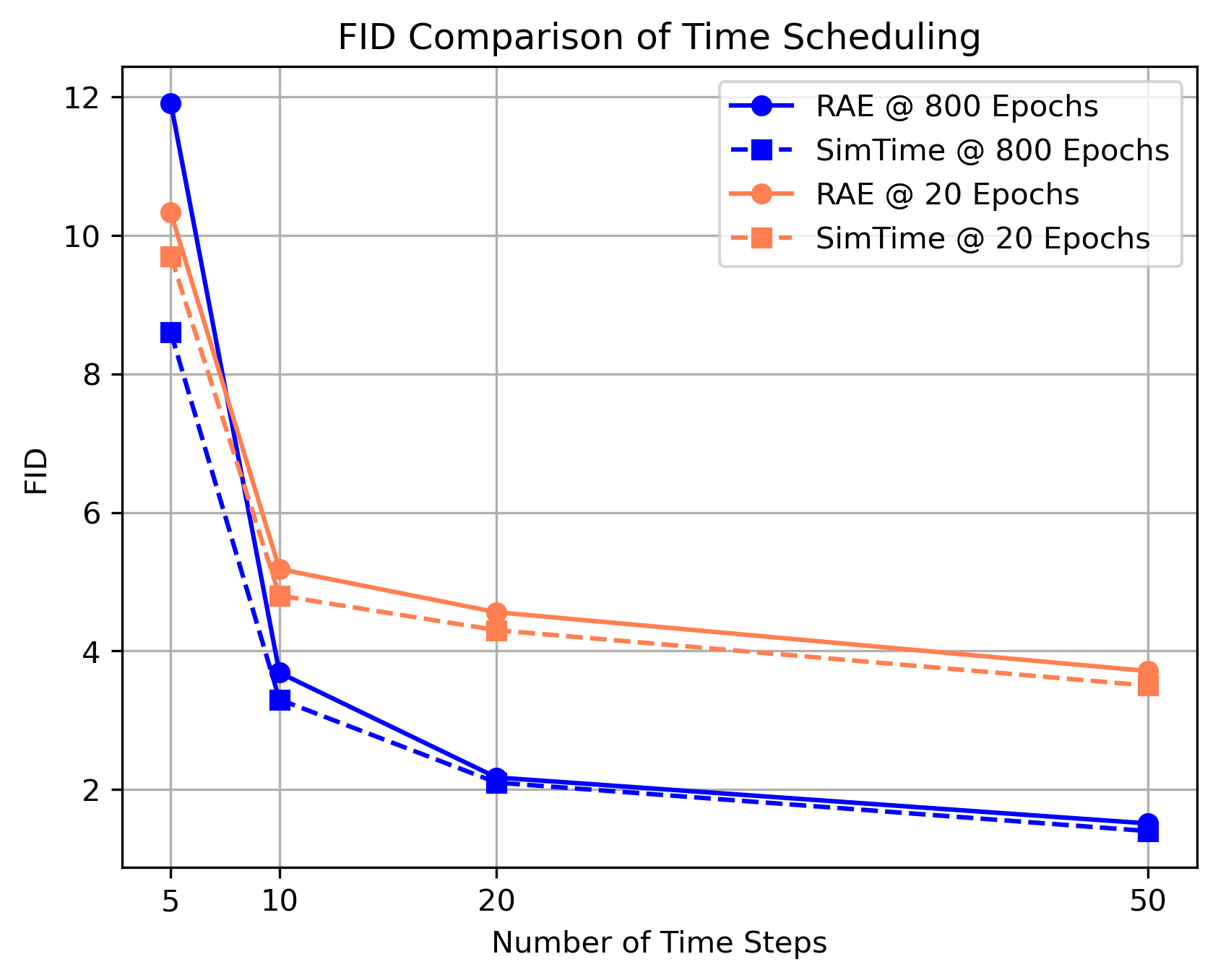}
			\caption{FID comparison of different time scheduling methods.}
			\label{fig:fid-timeschedule}
		\end{minipage}
		
	\end{figure}
	
	\section{Cosine Similarity Based Fine-Tuning Techniques}
	Motivated by recent advances in feature–alignment strategies for generative models
	\cite{tongImprovingGeneralizingFlowbased2024, liImmiscibleDiffusionAccelerating2024},
	we introduce a cosine-similarity–driven coupling scheme for fine-tuning.  
	This coupling selects data pairs by maximizing directional agreement in the latent or feature space, and the following proposition formalizes its optimality under an appropriate transport cost.  
	The proof is provided in Appendix~\ref{apdx:prof}.
	\begin{proposition}\label{prop:cos}
		Let $\mu$ and $\nu$ be probability measures on $\mathbb{R}^d$, and let 
		\[
		\Pi(\mu,\nu) := \{\gamma \text{ probability measure on } \mathbb{R}^d\times\mathbb{R}^d 
		: (\pi_X)_\#\gamma=\mu,\; (\pi_Y)_\#\gamma=\nu\}
		\]
		be the set of all couplings.
		
		Define the cosine similarity
		\[
		\phi(x,y) := \cos(x,y) = \frac{\langle x,y\rangle}{\|x\|\|y\|}
		\]
		and consider the cost function 
		\[
		c(x,y) := -\phi(x,y) = -\cos(x,y).
		\]
		
		Then any coupling $\gamma^\star \in \Pi(\mu,\nu)$ that maximizes the expected cosine similarity
		\[
		\gamma^\star \in \arg\max_{\gamma\in\Pi(\mu,\nu)} 
		\int_{\mathbb{R}^d\times\mathbb{R}^d} \phi(x,y)\, d\gamma(x,y)
		\]
		is also an optimal transport plan for the cost $c$, i.e.
		\[
		\gamma^\star \in \arg\min_{\gamma\in\Pi(\mu,\nu)} 
		\int_{\mathbb{R}^d\times\mathbb{R}^d} c(x,y)\, d\gamma(x,y).
		\]
		Conversely, any minimizer of the cost is a maximizer of the cosine similarity.
	\end{proposition}
	\begin{table}[H]
		\centering
		\caption{Comparison of generative models at 256$\times$256 resolution.}
		\small
		\begin{tabular}{lcccccc}
			\toprule
			\multirow{2}{*}{\textbf{Method}} &
			\multirow{2}{*}{\textbf{Epochs}} &
			\multirow{2}{*}{\textbf{\#Params}} &
			\multicolumn{2}{c}{\textbf{Generation @256 w/o guidance}} &
			\multicolumn{2}{c}{\textbf{Generation @256 w/ guidance}} \\
			\cmidrule(lr){4-5} \cmidrule(lr){6-7}
			& & & gFID$\downarrow$ & IS$\uparrow$ & gFID$\downarrow$ & IS$\uparrow$ \\
			\midrule
			\multicolumn{7}{l}{\textit{Autoregressive}} \\
			VAR \cite{tianVisualAutoregressiveModeling2024} & 350 & 2.0B & 1.92 & 323.1 & 1.73 & \textbf{350.2} \\
			MAR \cite{liAutoregressiveImageGeneration2024} & 800 & 943M & 2.35 & 227.8 & 1.55 & 303.7 \\
			\cite{renNextTokenNextXPrediction2025} & 800 & 1.1B & -- & -- & 1.24 & 301.6 \\
			\midrule
			\multicolumn{7}{l}{\textit{Latent Diffusion with VAE}} \\
			DiT \cite{peeblesScalableDiffusionModels2023} & 1400 & 675M & 9.62 & 121.5 & 2.27 & 278.2 \\
			MaskDiT \cite{zhengFastTrainingDiffusion2023} & 1600 & 675M & 5.69 & 177.9 & 2.28 & 276.6 \\
			SiT \cite{maSiTExploringFlow2024} & 1400 & 675M & 8.61 & 131.7 & 2.06 & 270.3 \\
			MDTv2 \cite{gaoMDTv2MaskedDiffusion2023} & 1080 & 675M & -- & -- & 1.58 & 314.7 \\
			VA-VAE \cite{yaoReconstructionVsGeneration2025} & 80 & 675M & 4.29 & -- & -- & -- \\
			& 800 & 675M & 2.17 & 205.6 & 1.35 & 295.3 \\
			REPA \cite{yuRepresentationAlignmentGeneration2025} & 80 & 675M & 7.90 & 122.6 & 1.29 & 306.3 \\
			& 800 & 675M & 5.78 & 135.8 & 1.29 & 306.3 \\
			DDT \cite{wangDDTDecoupledDiffusion2025} & 80 & 675M & 6.62 & 135.2 & 1.52 & 263.7 \\
			& 400 & 675M & 6.27 & 154.7 & 1.26 & 310.6 \\
			REPA-E \cite{lengREPAEUnlockingVAE2025} & 800 & 675M & 3.46 & 159.8 & 1.67 & 266.3 \\
			& 800 & 675M & 1.70 & 217.3 & 1.15 & 304.0 \\
			RAE \cite{zhengDiffusionTransformersRepresentation2025} & 20 & 839M & 3.71 & 198.7 & -- & -- \\
			& 80 & 839M & 2.16 & 214.8 & -- & -- \\
			& 800 & 839M & \textbf{1.51} & \textbf{242.9} & \textbf{1.13} & 262.6 \\
			\midrule
			\multicolumn{7}{l}{\textit{Fine Tuning (Ours)}} \\
			Cos & 20(+1) & 839M & 3.30 & 210.3 & -- & -- \\
			\bottomrule
		\end{tabular}
	\end{table}
	
	In the evaluation above, we fine-tune RAE for only one additional epoch beyond the 20-epoch checkpoint, yet obtain performance comparable to the 80-epoch baseline. This demonstrates that cosine-similarity–guided coupling provides an efficient fine-tuning signal with minimal computational overhead.
	\bibliography{reference}
	\bibliographystyle{plainnat}
	\appendix
	\section{Proof of Proposition~\ref{prop:cos}}\label{apdx:prof}
	\begin{proof}
		For any coupling $\gamma\in\Pi(\mu,\nu)$, define
		\[
		J(\gamma) := \int \phi(x,y)\, d\gamma(x,y),
		\qquad
		\mathcal{C}(\gamma) := \int c(x,y)\, d\gamma(x,y).
		\]
		By definition of $c$, we have
		\[
		\mathcal{C}(\gamma)
		= \int -\phi(x,y)\, d\gamma(x,y)
		= - J(\gamma).
		\]
		
		Let $\gamma_1,\gamma_2 \in \Pi(\mu,\nu)$ be arbitrary. Then
		\[
		J(\gamma_1) \ge J(\gamma_2)
		\quad\Longleftrightarrow\quad
		- J(\gamma_1) \le -J(\gamma_2)
		\quad\Longleftrightarrow\quad
		\mathcal{C}(\gamma_1) \le \mathcal{C}(\gamma_2).
		\]
		
		Thus, $\gamma^\star$ maximizes $J$ over $\Pi(\mu,\nu)$ if and only if it minimizes 
		$\mathcal{C}$ over $\Pi(\mu,\nu)$. Equivalently,
		\[
		\arg\max_{\gamma\in\Pi(\mu,\nu)} J(\gamma)
		=
		\arg\min_{\gamma\in\Pi(\mu,\nu)} \mathcal{C}(\gamma).
		\]
		This proves the claim.
	\end{proof}
	\section{Experimental Setup}
	The overall architecture follows that of \cite{zhengDiffusionTransformersRepresentation2025}. For the visual encoders, we employ DINOv2 with Registers \cite{darcetClusterPredictLatent2025}. We use the “large” configuration of DINOv2, which processes 224×224 images with a patch size of 14, yielding visual representations of dimension 768 over 256 tokens. For the diffusion backbone, we adopt a DiT model equipped with a shallow but wide DDT-head transformer module for denoising.
	
	We further incorporate several architectural enhancements, including SwiGLU feed-forward layers \cite{shazeerGLUVariantsImprove2020}, rotary positional embeddings \cite{suRoFormerEnhancedTransformer2021}, and RMSNorm \cite{zhangRootMeanSquare2019}.
	\begin{table}[t]
		\centering
		\caption{Hyperparameter settings across different model scales.}
		\label{tab:hyperparams_scales}
		\begin{tabular}{lccc}
			\toprule
			Backbone & DiT-B & DiT-L & DiT-XL \\
			\midrule
			\multicolumn{4}{l}{\textbf{Architecture}} \\
			\midrule
			\#Params      & 193M            & 617M            & 844M            \\
			Input         & $16\times16\times768$ & $16\times16\times768$ & $16\times16\times768$ \\
			Layers (DiT)        & 12              & 24              & 28              \\
			Layers (DDT Head)      & 2              & 2              & 2              \\
			Hidden dim. (DiT)   & 768             & 1024         & 1152         \\
			Hidden dim. (DDT Head)   & 2048             & 2048         & 2048   \\
			Num. heads (DiT)  & 12              & 16              & 16              \\
			Num. heads (DDT Head)   & 16              & 16              & 16              \\
			\midrule
			\multicolumn{4}{l}{\textbf{Optimization}} \\
			\midrule
			Batch size       & 256           & 256           & 256           \\
			Optimizer        & AdamW         & AdamW         & AdamW         \\
			lr               & 0.0001        & 0.0001        & 0.0001        \\
			$(\beta_1,\beta_2)$ & (0.9, 0.95) & (0.9, 0.95) & (0.9, 0.95) \\
			\midrule
			\multicolumn{4}{l}{\textbf{Interpolants}} \\
			\midrule
			$\alpha_t$            & $1-t$          & $1-t$          & $1-t$          \\
			$\sigma_t$            & $t$            & $t$            & $t$            \\
			$w_t$                 & $\sigma_t$     & $\sigma_t$     & $\sigma_t$     \\
			Training objective    & v-prediction   & v-prediction   & v-prediction   \\
			Sampler               & Euler & Euler & Euler \\
			Sampling steps        & 50            & 50            & 50            \\
			\bottomrule
		\end{tabular}
	\end{table}
	
	\section{Evaluation Details}
	We follow the evaluation protocol of \cite{liAutoregressiveImageGeneration2024} and employ the same reference statistics provided in their official implementation.\footnote{\url{https://github.com/LTH14/mar/tree/main}}
	For each experiment, we generate 50k samples using the SDE Euler–Maruyama sampler with 50 discretization steps.
	
	We report two standard metrics for generative image quality, both computed using features extracted from the Inception-V3 network. Following common practice, all features are obtained from the Inception-V3 classifier pre-trained on ImageNet \cite{szegedyRethinkingInceptionArchitecture2015}.
	
	\begin{itemize}
		\item \textbf{Inception Score (IS)}~\cite{salimans2016improved}: Evaluates both sample fidelity and diversity by measuring the KL divergence between the conditional class distribution and the marginal class distribution.
		\item \textbf{Fréchet Inception Distance (FID)}~\cite{heusel2017gans}: Quantifies the discrepancy between generated and real image distributions by computing the Fréchet distance between their Inception-V3 feature statistics.
	\end{itemize}
	
	\section{Visual Results}
	\begin{figure}
		\includegraphics[scale=0.4]{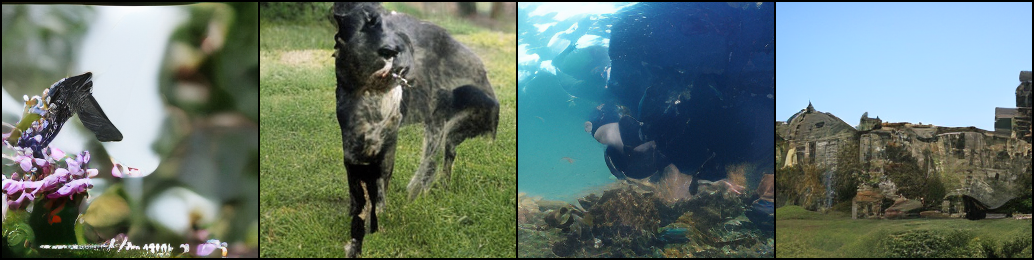}
		\includegraphics[scale=0.4]{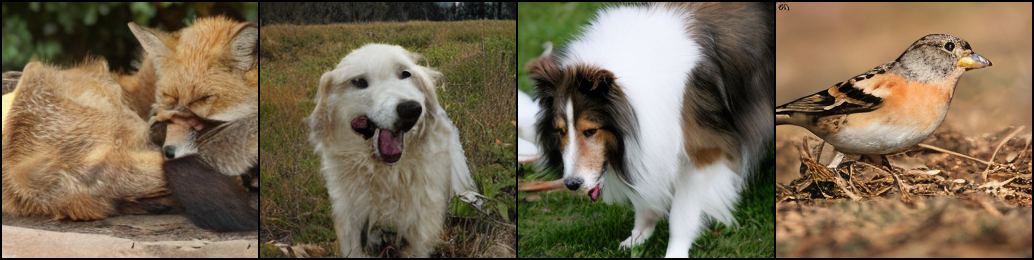}
		\caption{Visual Results for Fine-Tuned Model}
	\end{figure}
\end{document}